# Analysis of Growing Tumor on the Flow Velocity of Cerebrospinal Fluid in Human Brain Using Computational Modeling and Fluid-Structure Interaction

Muhammad Uzair-Ul-Haq[1], Ali Ahmed[1], Zartasha Mustansar[2,*], Arslan Shaukat[3], Lee Margetts[4], Asim Waris[1] and Faizan Nadeem[5]

[1] Department of Biomedical Engineering & Sciences, School of Mechanical & Manufacturing Engineering, National University of Sciences and Technology, Islamabad 44000, Pakistan;
email: ulhaquzair@yahoo.com; talpur54@outlook.com
[2] Department of Computational Engineering, Research center of Modelling and Simulation, National University of Sciences and Technology, Islamabad 44000, Pakistan; email: zmustansar@rcms.nust.edu.pk
[3] Department of Computer & Software Engineering, College of Electrical & Mechanical Engineering, National University of Sciences and Technology, Islamabad 44000, Pakistan; email: arslanshaukat@ceme.nust.edu.pk
[4] Department of Mechanical, Aerospace and Civil Engineering, University of Manchester UK; email: lee.margetts@manchester.ac.uk
[5] Nuclear Medicine specialist associated with INMOL cancer hospital, Lahore, Pakistan; email: dr_faizi@hotmail.com
* Correspondence: zmustansar@rcms.nust.edu.pk.

**Abstract**: Cerebrospinal fluid (CSF) plays a pivotal role in normal functioning of Brain. Intracranial compartments such as blood, brain and CSF are incompressible in nature. Therefore, if a volume imbalance in one of the aforenoted compartments is observed, the other reaches out to maintain net change to zero. Whereas, CSF has higher compliance over long term. However, if the CSF flow is obstructed in the ventricles, this compliance may get exhausted early. Brain tumor on the other hand poses a similar challenge towards destabilization of CSF flow by compressing any section of ventricles thereby ensuing obstruction. To avoid invasive procedures to study effects of tumor on CSF flow, numerical-based methods such as Finite element modeling (FEM) are used which provide excellent description of underlying pathological interaction. A 3D fluid-structure interaction (FSI) model is developed to study the effect of tumor growth on the flow of cerebrospinal fluid in ventricle system. The FSI model encapsulates all the physiological parameters which may be necessary in analyzing intraventricular CSF flow behavior. Findings of the model show that brain tumor affects CSF flow parameters by deforming the walls of ventricles in this case accompanied by a mean rise of 74.23% in CSF flow velocity and considerable deformation on the walls of ventricles.

*Keywords: Intracranial pressure, cerebral edema, non-invasive treatment, brain tumor, fluid-structure interaction.*

# 1. Introduction

*1.1. Importance and Rational of Study*

Biomechanical modeling of Cerebrospinal fluid (CSF) flow is one of the important factors which merits investigation in the area of neuropathology. CSF plays a pivotal role in brain physiology. Continued forces of brain tumor on the brain tissue can affect the normal working of brain. In a current scenario where a tumor is present in the close vicinity of brain ventricles, where CSF flows normally, the increasing tumor mass over a period of time can deform the brain ventricles hence disturbing the flow of CSF and the homeostasis of brain. In the aforementioned case, brain tumor induces stenosis of cerebral aqueduct causing obstructive hydrocephalus condition. Clinical symptoms of stenosis of cerebral aqueduct, for instance include seizures, gait imbalance, visual problems [3]. It is also important to note that in patients having complaint of elevated ICP problem, one of the leading causes is the obstruction to the flow of CSF [4].

Therefore it is imperative, and necessary, to describe by means of a numerical biomechanical model the underlying interactions of brain tumor and CSF flow with the ventricular walls which is critical in developing anatomically realistic model of CSF flow mechanics. The focus of this investigation would thus revolve around developing a complete biomechanical model of CSF in brain ventricles and tumor growth over a period of time which may further aid clinicians and practitioners in taking informed decisions before undertaking any surgical intervention.

*1.2. CSF Biomechanics*

CSF is a colorless incompressible and Newtonian fluid which is found both in the human ventricular system (HVS) and in the Subarachnoid Space (SAS) [10]. It is produced in the Choroid Plexus, found inside on the walls of Ventricles. Ventricles are cavities or sac like structures and function as containers for the CSF. There are in total two lateral ventricles, one in each cerebral hemisphere [10]. Lateral ventricles are connected to each other via Intraventricular Foramina which is the junction through which CSF flows into the Third Ventricle; it then extends into fourth ventricle via Aqueduct of Sylvius and thereafter exits from two sets of apertures viz. Foramen of Luschka and Foramen of Magendie. Penultimately, the CSF flows and circulates both in the Spinal cord and in the SAS. At nominal conditions, the volumetric production rate of CSF approximates to 500 ml/day [10]. The flow direction of CSF inside cranium has two distinct stages within a cardiac cycle. During systole, brain is compressed due to the pressure of blood forcing in through cerebral vasculature thereby making flow through lateral and third ventricle in craniocaudal direction [11]. While during diastole, the flow through the third ventricle and cerebral aqueduct reverses due to contraction of blood vessels moves. Hence contraction and expansion of blood vessels (or cerebral vasculature) creates the necessary motion of CSF down the cerebral aqueduct and fourth ventricle all the way till the spinal canal. This to and fro exercise also causes pulsatility in CSF.

*1.3. Biomechanical Modeling of CSF—Research Gap*

In recent decades, numerical methods along with computer simulations are commonly being used to discretize real life scenarios such as blood flow in arteries, cardiac muscle modeling. Computational fluid dynamics (CFD) has been widely used to model CSF flow. Analysis of CSF behavior and properties such as velocity profile and pressure drop through ventricles under influence of any pathology allow the clinicians to differentiate between normal and abnormal cases. CSF modeling has been of great interest by researchers from across the globe. For instance, one of the proposed models [12] consists of a simple cylinder model to study the flow of CSF through of cerebral aqueduct. Afore cited study suggests that pressure difference of at least 1.1 Pa is required to drive the CSF from cerebral aqueduct. Fin et al [13] also proposed two models, one cylindrical rigid wall model and other an elastic wall model segmented from MRI data to study the flow of CSF in Aqueduct of Sylvius (AS). The spatial domain was discretized using Immersed Boundary method (IBM) [14]. Pulsatile inlet velocity boundary condition is used, and their results show pressure drop of 1.02 Pa and velocity of 30.20 mm/s for cylindrical model; whereas for elastic wall model, the pressure drop of 2.91 Pa and velocity of 64.65 mm/s was recorded. In another case study Linninger et al. [15] proposed a 2D rigid model of HVS where the pressure drop of less than 2 Pa and velocity of 7.3 mm/s is reported. In yet another case study, Kurtcuoglu et al. [16] constructed a simplified geometry of HVS to model the mechanical behavior of CSF. To drive the flow of CSF, oscillating wall of third ventricle was used whereas at the outlet, zero-pressure boundary condition was applied. The study deduced that pressure increases in lateral ventricles due to aqueduct stenosis but no significant phase difference between wall motion and pressure was observed. Moreover, the CSF movement is influenced by the pathway of CSF through CNS and brain movement. In addition, the pulsatile nature of CSF is also due to the elastic nature of the Ventricular space, and due to the changes in systolic and diastolic pressure in each cardiac cycle [17]. One of the extensive studies done on the subject matter in literature is reported by Howden et al. [18] wherein author used anatomically realistic ventricles which were 3D Ventricles segmented from MRI data with 1 mm voxel size using MIMICS. A commercial code of Fluent was used for CFD analysis. CSF is treated as incompressible Newtonian fluid. Pulsatile velocity inlet and zero-gauge pressure at outlet boundary conditions were used. The maximum velocity and Reynolds number was found to be 11.38 mm/s and 15 in cerebral aqueduct respectively. However, in this study the walls of ventricles are taken as rigid and effects of Fluid-structure interaction are neglected.

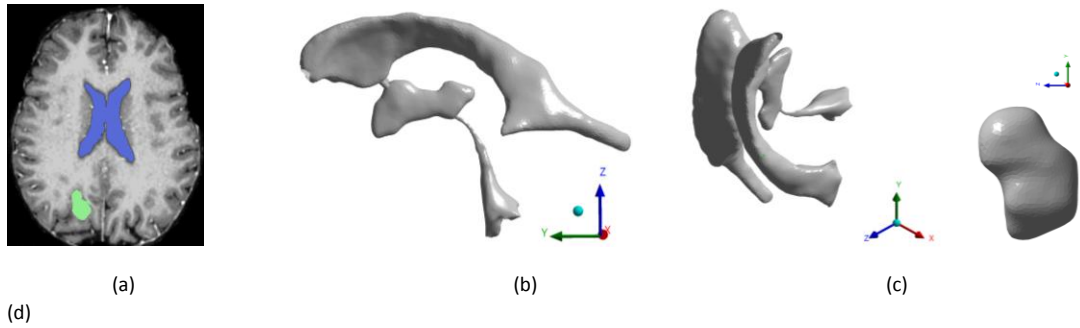

**Figure 1**. Representation of ventricles; (a) shows Axial section of MRI image of human brain showing ventricles and tumor; (b) and (c) show 3D reconstructed ventricular geometry in lateral and isometric view; (d) show 3D reconstructed tumor geometry.

In past CSF modeling has been limited to a simple flow analysis by taking ventricular walls as rigid, ignoring completely the influence of deformable nature of ventricles, thus hindering in towards the reaching towards a complete understanding of flow dynamics and tissue interaction mechanics. This modeling approach is also problematic due to the fact that by considering the walls as rigid, effect of the external forces such as that of brain tumor cannot be included into broader framework as of how tumor affects the CSF flow parameters within the brain ventricles. Therefore, a course correction in approach for modeling such situations is needed. Hence, we present in this paper a full 3-Dimensional biomechanical model of CSF flow which accurately captures both the flexible nature of ventricles and tumor growth effect on CSF velocity and pressure profile using Finite element and FSI modeling.

## 2. Materials and Methods

### 2.1. 3D Geometry and MRI Dataset

In this study, three sets of MRI scans are used. One datapoint is based on a real 3D T1-contrast enhanced MRI scan of a tumor patient (Datapoint 1), second and third datapoints (Datapoint 2 and 3) consists of two set of MRI scans taken from BRATS 2018 dataset [27]. MRI data in the former case is obtained with a 1.5 Tesla MRI machine and has an image resolution of 512 mm x 512 mm x 288 mm. While the latter data points (BRATS 2018 dataset) have image resolution of 240 mm x 240 mm x 155 mm. Two different data sources are used to validate the results on different specification scales and subsequently showing the performance of our proposed model.

Ventricular and tumor geometries are segmented using 3D Slicer and Mimics Innovation Suite. Semi-automatic method of segmentation is used to segment the geometries and special care is given to preserve image features which represent sharp curvatures or textural intricacies. The segmented geometries are compared with the ground truth provided by a Radiologist to validate the segmentation process. The calculated segmentation volumes and ground truths are provided in Table 2. Figure 1 (a) and (b) show axial scans highlighting ventricle and tumor region. For smoothing of Computer Diagnosis (CAD) models, Laplacian smoothing filter is applied, and 3D CAD models are exported into ANSYS (analysis-based software for Finite element operations) for further analysis. 3D volume is converted into surface mesh of triangular shapes using Delaunay triangulation scheme [37]. Figure 1 (c) and (d) shows reconstructed ventricle and tumor geometry.

*2.2. Meshing and Mesh Independence test*

Mesh independence study is performed to find the ideal grid size. As Fluid-Structure Interface technique is used therefore performing mesh independence test on fluid and structure domain separately gives better indication of the converged solution. In solid domain, on each mesh size a static structural solution was applied out under a specified pressure and response thereof in terms of peak stresses was noted. By performing the mesh independence test it was found that the grid size of 196576 elements with minimum element size of 0.3 mm is fine enough to capture all the physics in the simulation as deformation and stresses especially around the cerebral aqueduct and third ventricle. Average stresses below the threshold grid size almost remain constant and no further variation is found in the convergence of stress as shown in Figure 2 (a). Similarly, in the case of fluid domain, tetrahedral elements (SOLID 187) are used for meshing since the geometry is complex and tetrahedral elements can better capture complex and intricate features. Mesh independence study was carried out to find out an ideal grid size. Figure 2 (b) shows that the 1144379 numbers of elements with minimum element size of 0.5 mm is fine enough to capture all the important information regarding peak velocity changes in the simulation.

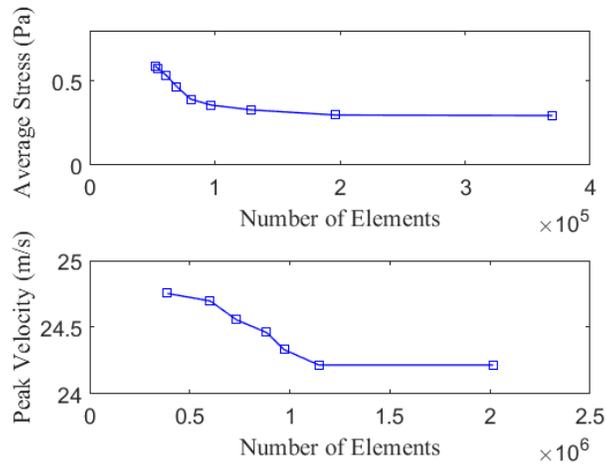

**Figure 2**. Grid independence study, (a) is for the structural domain, and (b) for the fluid domain.

*2.3. Properties of Ventricle Body and CSF fluid*

The ventricles are modeled as elastic [19] and also have viscous damping effect as provided by [R]. The elastic nature of the material is modeled using simple stress-strain constitute law where stress is directly proportional to strain:

$$\sigma = E\varepsilon \tag{1}$$

*where stress is directly proportional to strain; here E is the modulus of elasticity, σ is a stress tensor, ε is a strain tensor*. The viscous property is applied using Prony Series shear relaxation given by:

$$G(t) = G_\infty + (G_0 - G_\infty)e^{\frac{-t}{\tau}} \qquad (2)$$

*Where $G_\infty$ is the long-term shear modulus, $G_0$ is the instantaneous shear modulus and τ is the time decay constant.* The CSF is modeled as Newtonian and incompressible fluid [28] having constant density equal to 1000 kg/m3. The flow profile is laminar in nature and therefore viscous laminar model in fluent is used. In throughout fluid domain, CSF flow behaves as a creeping flow [28] therefore viscous effects dominate, and inertial effects do not affect much. Viscosity of CSF is taken as 1.003 mPa.s [19]. Table 1 enlists the material properties used in this study.

**Table 1** Material properties

| Boundary Conditions | Value | Reference |
|---|---|---|
| CSF Density | 1000 kg/m³ | Masoumi et al [19] |
| CSF viscosity | 0.001003 Ns/m² | Redzic et al [41] |
| Bulk mass flow | 6.25x10⁻⁶ kg/s | Wright et al [10] |
| Pressure outlet | Zero Pascal | L. Howden et al [18] |
| Modulus of Elasticity | 30 KPa | Masoumi et al [19] |
| Viscous Property | G(t)= 0.00101+(0.101-0.00101)e^-t/100 | Zhang et al [48] |
| Ventricles Density | 1000 kg/m³ | Masoumi et al [19] |
| Poisson's ratio | 0.49 | Masoumi et al [19] |

*2.4. Boundary Conditions*

As shown in Figure 3, the model consists of two inlets and two outlets. Inlets are defined in lateral ventricles with 3mm holes. Outlets are defined in foramen of Magendie (lateral aperture) and Luschka (median aperture). Mass flow inlet boundary condition with bulk production of 500 ml/day or 6.25x10-6 kg/s [19] is used to mimic the real-life scenario of CSF production. While at the outlet, zero-gauge pressure boundary condition is used to make sure that the outlet remains at a fixed static pressure. The inlets are specified at the lateral ventricles because since the overall flow profile inside ventricle is of creeping flow.

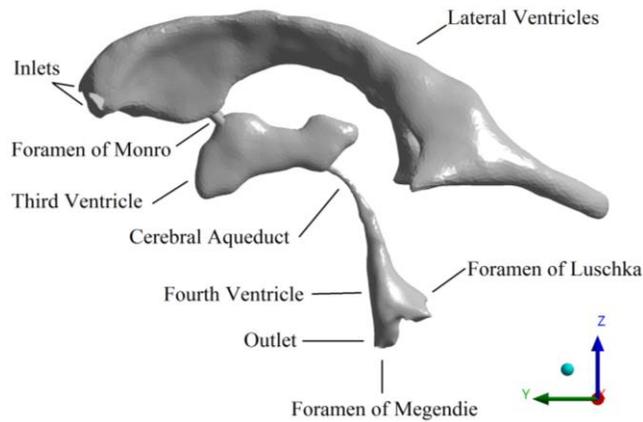

**Figure 3**. Boundary conditions

*2.5. CSF Flow Pulsatility*

As noted above, CSF flow has to and fro motion due to artifacts of blood pressure in a cardiac cycle. Hence physiologically, CSF production has a constant bulk production detailed above, but it has also a pulsatile component due to cardiac-induced systole-diastole movement. The pulsatile component is ideally a linear combination of sinusoidal harmonics and the combined form is given as below.

$$V(t) = B + X\sin(wt + \alpha) + Y\sin(2wt + \beta) \qquad (3)$$

*where B is the bulk production which equals to 0.3 ml/day, X is equal to 0.21 ml/min, and Y is the part of second harmonic which equals to 0.05ml/min* [29].

*2.6. Modeling of Tumor Forces*

In human pathology, brain tumor is characterized by the abnormal cells which gradually grows over time and affects the physiology of brain. To study the effect of tumor growth on the human ventricles system and flow of cerebrospinal fluid, the tumor growth is predicted using Gompertz mathematical model which is one of the widely used models to extrapolate tumor cells based on the initial tumor volume [30]. The mathematical model has the solution in the form as:

$$V_f(t) = V_o e^{\frac{\alpha}{\beta}(1-e^{-\beta t})} \qquad (4)$$

*Where $V_o$ is the initial tumor volume, $V_f$ is the final tumor volume, α and β are the initial proliferation rate and the exponential decay constants and t is the time in days. The value of α and β are taken as 0.279 and 0.1470 respectively* [32]. While the efficacy of these models on humans is still an ongoing effort and needs to be gauged objectively because of the inadequacy of in-vivo measurements of tumor size over time, various models such as those aforementioned, however have mostly relied on animal studies such as those done on mice to find the constant values.

With regards to tumor forces, the tumor content is not fully a solid core [33], it also contains hydrostatic forces of fluid. However, in a malignant case the tumor mass mostly consists of solid core therefore, for simplicity the action of fluid forces is ignored, and only solid core is considered (which is also the extracellular matrix of the tumor). Body forces is one parameter which can address the concern of finding the forces of solid tumor core. Body forces can be thought of as forces acting on entire volume of body such as those forces due to the gravity. Forces due to gravity are basically the weight of the body; whereas tumor forces are applied on the walls of third ventricle and can be given by equation (9) as:

$$F = \int_V \rho g dV \qquad (5)$$

*where F is the body force, ρ is the density of the tumor mass and dV is the volume differential element. In the above equation ρ is taken as 1040 kg/m3* [34] *and g is taken as 9.81 m/s2.*

The predicted tumor forces for all the datapoints are given in Table 2. The tumor forces are applied on the walls of third ventricle. The aim of this study is to model how tumor core as a whole interacts with the walls of ventricles and effects the cerebrospinal fluid. Such a case can be visualized in scenarios wherein there is an intraventricular tumor or a tumor in very close proximity such as in brain stem. In the current scenario tumor applies relatively maximum force on the walls of ventricles. This is a specific class of cases and we aim to model them by providing a fundamental framework to study them using

| Specimen | Predicted Volume (mm³) | Predicted Force (N) |
|---|---|---|
| Datapoint 1 | 6590 | 0.0672 |
| Datapoint 2 | 3604 | 0.0368 |
| Datapoint 3 | 8592 | 0.0877 |

Table 2 Predicted Tumor Volumes and Forces.

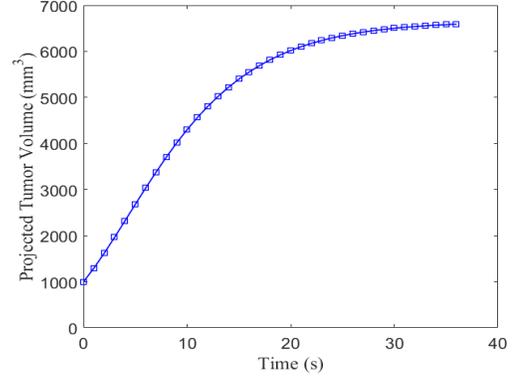

Figure 4. Predicted tumor volume Grown over 35 days

FSI approach. The proposed model can thereafter be further extended to include brain parenchyma interaction if the tumor's spatial location is not in the vicinity of ventricles and general brain deformation aspects need to be catered. Figure 4. shows projected tumor growth for datapoint 1 (here time is in days and is only used for growth purposes).

### 2.7. Mathematical Method

The biomechanical model of CSF flow in ventricles is modeled in ANSYS using fully coupled 3D fluid structure interaction technique. The fluid flow is solved using Navier-Stokes equations. Mass and momentum conservation equations are given as below [35]:

$$\nabla \cdot (\rho \cdot \vec{v}) = 0 \qquad (6)$$

$$P\frac{\partial \vec{v}}{\partial t} + \rho(\vec{v} \cdot \nabla)\vec{v} = -\nabla P + \tau\rho + \mu\nabla^2 \vec{v} \qquad (7)$$

where $\rho$ is the density of fluid, v is the velocity of the fluid, $\nabla$ is the gradient operator, $\rho\frac{\partial \vec{v}}{\partial t}$ is the local acceleration of fluid particles, $\rho(\vec{v}.\nabla)\vec{v}$ is the convective acceleration, $\nabla P$ is the pressure gradient, $\tau\rho$ are body forces, $\mu\nabla^2\vec{v}$ is the viscous term which resists the motion of the fluid particles. . Solving above equations (6) and (7) numerically requires spatial discretization and pressure-velocity coupling schemes which can be used to interpolate pressures at the faces of a control volume. Pressure-based solver is used with PISO (Pressure-Implicit with Splitting of Operators) coupling is used with momentum discretization scheme used as second order upwind; whereas for gradient discretization least-square cell-based method is used. Implicit formulation is used in discretizing equations and equations are linearized and solved iteratively. This results into N equations for each cell. An implicit equation solver, based on Incomplete-Lower-Upper (ILU) factorization scheme, together-with Algebraic Multigrid method (AMG), is used to speed up the convergence of the solution.

To capture the physics of stress and deformation due to fluid flow, Transient Structural is used where solid domain is numerically solved. Newmark Integration scheme is used to solve equations because the model generated in our scheme has non-linearities in it, therefore implicit solver is needed to reach towards a solution in an iterative manner. To judge the convergence of solution Newton-Raphson method is used which is given as below:

$$x_{n+1} = x_n - \frac{f(x_n)}{f'(x_n)} \qquad (8)$$

Equation (8) is used to converge forces, moments and displacements at each time-step, which works by the principle that the energy added due to the external loads must eventually balance the energy induced by the reaction forces. Due to deformable boundary the elements can get distorted and thus losing their ideal shape, therefore dynamic meshing is needed so that the moving or deforming cells can be re-meshed correctly thereby preserving shape quality. Remeshing is used to keep the mesh parameters such as skewness, aspect ratio and orthogonal quality in acceptable quality range to provide better results. Diffusion-based smoothing is used to smooth deforming cells.

$$\nabla \cdot (\gamma \nabla \vec{u}) = 0 \tag{9}$$

$$\gamma = 1/d^a \tag{10}$$

*where d is the normalized boundary distance, a is the diffusion parameter and $\vec{u}$ is the mesh displacement velocity.* To ensure true convergence, timestep needs to be same in both the solvers so that data exchange occurs at same instance. A possible scenario for calculating a reasonable timestep size is by considering the cardiac-induced pulsations which would amount to taking smallest timestep which is calculated using Δt=1/20f, where Δt is the time-step size, and f is the frequency of heart beats. For a normal person, heartbeat ranges from 60-100 bpm. Frequency of 60 bpm equals 1 Hz and using aforementioned formula Δt equals to 0.05s. In this study we have taken a time-step of 0.01 seconds which is smaller than the above calculated threshold to capture small changes in the simulation. Lastly, relative convergence in the structural domain (Transient Structural), fluid domain (Fluent) and in the system coupling is achieved so that results correspond to true convergence. Convergence in data exchanges is achieved and is set to root-mean-square change of 0.0005.

## 3. Results

The numerical model is simulated under two broad cases: first case relates to the normal scenario where no tumor is present and normal CSF flow in ventricle is studied. The second case relates to the interaction of brain tumor to the ventricular body. In the second case the effect of tumor growth on the functioning of human ventricles is studied. In this case two sets of interactions are observed: first the effect CSF flow on the walls of ventricles in the presence of tumor forces are captured and secondly, the resulting deformation of ventricular walls and their effect on the flow of CSF in the presence of tumor. The 3D anatomically realistic geometrical model of ventricles and brain tumor so segmented is validated in terms of the ventricular and tumor volumes against the ground truths obtained through radiologist which are supplied in Table 3.

**Table 3** Ventricular and Initial Tumor Volume Validation Against Ground Truths

| Specimen | CSF Ventricular Volume | | | Initial Tumor Volume | | |
|---|---|---|---|---|---|---|
| | Segmented Volume (mm³) | Ground truth Volume (mm³) | Percentage error (%) | Segmented Volume (mm³) | Ground truth Volume (mm³) | Percentage error (%) |
| Datapoint 1 | 25320 | 24320 | 4.11% | 997.29 | 959.51 | 3.93% |
| Datapoint 2 | 26940 | 25310 | 6.44% | 545.36 | 513.41 | 6.22 % |
| Datapoint 3 | 26760 | 25080 | 6.69% | 1300.14 | 1247.30 | 4.23 % |

*3.1. CSF Flow Velocity Due to Tumor Interaction*

The results obtained from the analysis are presented in figure 5 (a) and (b) Peak CSF velocity occurs in the aqueduct of Sylvius due to narrowing of pathway. Findings of present study show that CSF velocity values in the aqueduct rise to a greater value and remain almost same in the foramen of Monro and lateral ventricles. As clear from figure 5 (a), maximum velocity of 22.84 mm/s occurs in aqueduct of Sylvius when no tumor forces are present; while for the case of tumor it rises to 39.40 mm/s as shown in figure 5 (b) indicating a rise of 72.50 %. This change in velocity is one of the significant parameters which qualitatively defines the CSF flow field, especially in the cerebral aqueduct where most changes take place. Table V shows the results obtained from Datapoint 1,2 and 3. Mean rise in velocity for cases of tumor is found out to be 74.23%.

Reynold's number is an important parameter which signifies the transition of flow from laminar to turbulent. It is defined as the ratio of fluid inertial forces to the viscous forces, and is a metric used to differentiate whether CSF, in localized region such as aqueduct, remains laminar or turbulent. It is given as:

$$Re = \frac{\rho v D}{\mu} \qquad (11)$$

*where D is the characteristic length which is taken to be the diameter of the aqueduct, μ is the dynamic viscosity of fluid, ρ is the fluid density and v is the velocity of the CSF in aqueduct.* Based on the above findings related to CSF velocity, maximum Reynold's number across normal cases is 60.47, and for tumor case it rises to 107.91. Overall flow remains in the laminar regime, but in localized sections, velocity gradients are abrupt.

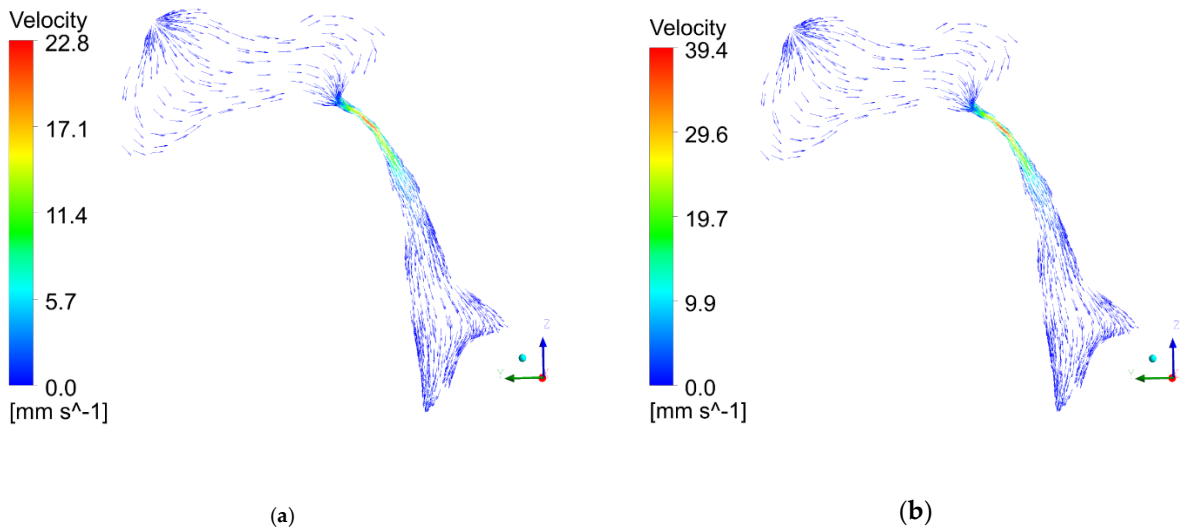

(**a**) (**b**)

Figure 5 Velocity vector plots in third and fourth Ventricle; (a) represents the case of no tumor, and (b) represents the case with tumor.

Transition of CSF from creeping flow towards a relatively high local Reynold number flow is evident for both cases. First, CSF flows from the lateral ventricles and passes into the foramen of Monro. Here, the velocity changes are insignificant and almost remain a creeping flow. It also remains the same for the tumor case as well in lateral ventricles. Recirculation of flow of small magnitude also appear in the third ventricle due to its geometry. When CSF enters into the cerebral aqueduct region, velocity starts to rise and reaches up to 22.84 mm/s for normal case, and for tumor case it reaches up to 39.4 mm/s clearly showing that the flow no longer, at least in the localized region of the cerebral aqueduct, remain strictly under the characterization of a creeping flow. Post aqueduct, velocity in the fourth ventricle in both cases and drops and flow goes back towards a relatively creeping flow regime. Reynold's number throughout the ventricular domain remains under 0.015 (considering velocity at outlet) showing that the flow by and large remains a creeping flow. Major oscillations in flow field are observed in the cerebral aqueduct. According to the phase-contrast MRI studies (in pre-shunting cases) [38] the velocity in normal cases remains under 39±13 mm/s whereas for the cases of stenosis of aqueduct due to mesencephalic tumor or non-communicating hydrocephalus it rises to 69±19 mm/s. This well correlates within the provided range and confirms the FSI-based biomechanical model presented in this paper. Reynold's number across all datapoints are enlisted in Table 4.

Table 4 CSF Velocity, Reynold's and Womersley number computed in Aqueduct of Sylvius.

| Specimen | CSF velocity | | | Reynold's number | | Womersley number |
|---|---|---|---|---|---|---|
| | Without tumor (mm/s) | With tumor (mm/s) | Percentage change (%) | Without Tumor (mm/s) | With Tumor (mm/s) | |
| Datapoint 1 | 22.84 | 39.40 | 72.50 | 41.0 | 70.7 | 5 |
| Datapoint 2 | 24.31 | 41.75 | 71.74 | 52.1 | 89.5 | 5.38 |
| Datapoint 3 | 28.97 | 51.70 | 78.46 | 60.47 | 107.91 | 5.75 |

*3.2. CSF Flow pulsatility under Cardiac Cycle*

To calculate pulsatile nature of flow field, Womersley number denoted by $\alpha$ is used. Womersley number is a dimensionless number and a metric used to signify the pulsatile nature of the biofluid in relation to the flow viscous effects [38]. It is given as below:

$$\alpha = L\sqrt{\frac{2\pi\rho}{T\mu}} \qquad (12)$$

where L is the characteristic length, T is the time period of the heartbeat oscillations, ρ is the fluid density and μ is the dynamic viscosity of the fluid. In case of ventricular system, characteristic length L is taken to be the diameter of aqueduct, $\alpha$ is calculated for all the datapoints as shown in the Table 4. For values of $\alpha \geq 1$ the flow field is sinusoidal in nature. Across all the datapoints, calculated Womersley number is greater than one, suggesting that the flow is pulsatile in nature. As stated earlier, pulsatile mass flow inlet boundary condition is used using equation (3). To visualize the effects of cardiac cycle and systole-diastole effects, the entire simulation is executed for 4 seconds. Figure 6 shows the output of pulsatile component of mass flow rate, peak

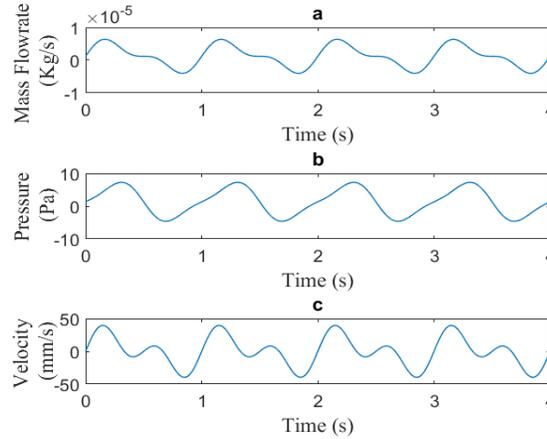

Figure 6. (a) pulsatile mass flowrate at inlets, (b) peak pressures in the lateral ventricles, and (c) CSF velocity in cerebral aqueduct.

pressures in the lateral ventricles, and CSF velocity in the cerebral aqueduct.

Human heart acts as a pumping machine and through entire cardiovascular system it supplies blood to all parts of body. Carotid arteries are the supply lines for blood towards brain. Since heartbeat is continuous and is usually rhythmic and periodic, pulsations in blood are formed which also affect the CSF movement in the ventricles and subarachnoid space. Every beat induces a to and fro motion of CSF in the ventricles. This motion is usually modeled by taking it as combination of sinusoidal waves whose frequency is equal to the frequency of heartbeat. Peak systole represents the maximum pressures the heart applies to supply blood to various parts of the body, and peak diastole is point where heart relaxes and refills the blood in its chambers. As visible from the Figure. 6, the maximum velocity occurs at approximately 17% (where frequency of beat is 1 Hz therefore, in window of 0-1 seconds, 17% corresponds to time of 0.17 seconds) of systole and the minimum velocity occurs at about 84% of the diastole.

*3.3. CSF Pressure Distributions*

Pressure changes in the ventricular wall represents the intracranial pressures. Pressure variations in large cavities such as lateral ventricles are negligible and found to be spatially uniform. Figure 7 (a) and (b) shows the pressure field for normal and tumor cases. Maximum pressures under the normal case are 4.35 Pa in the lateral ventricles. Pressure drop of 0.1 Pa is found in the third ventricle, while across the cerebral aqueduct it is found out to be 2.88 Pa. While for the tumor case, maximum pressure of 7.1 Pa is found in the lateral ventricles, with a pressure drop of 0.25 Pa in the third ventricle. Pressure drop in the cerebral aqueduct is found out to be equal to 4.74 Pa. A large pressure drop in the cerebral aqueduct shows a higher velocity rise. Comparing pressures found for normal case with some previous FSI-based 3D models, pressures on large cavities remain within the range of 4 Pa which validates the present results obtained [20].

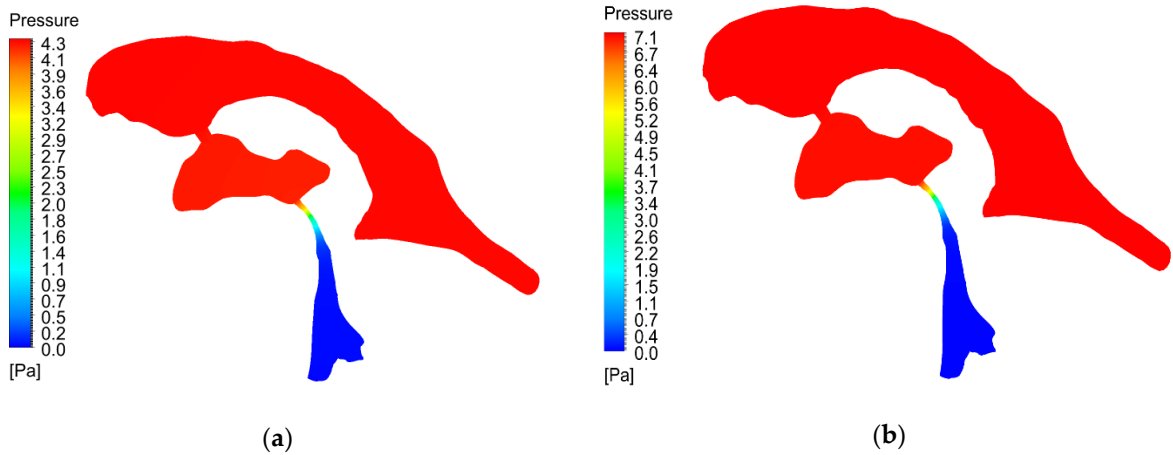

(a)                  (b)

Figure 7. Pressure on ventricle wall. (a) shows pressures without tumor; and (b) shows pressure under tumor forces

*3.4. CSF Backflow Visualization*

CSF flows back and forth in the ventricles in a cardiac cycle. During systolic pressure, lateral and third ventricles are compressed thereby forcing the CSF movement [11] towards the spinal canal and in subarachnoid space; at this instance velocity is on the higher side even in the cerebral aqueduct and in the fourth ventricle. Figure. 8 shows velocity vector plots for two cases: case A refers to normal case where no tumor is present and case B refers to a tumor case. At time steps 0.25 and 0.5 for case A, one can see that the velocity peaks at 0.25s, and post 0.5s, the diastolic phase starts, and velocity starts to reverse; timesteps of 0.75s and 1s show that velocity starts to dip and then peaks again respectively. The same is repeated for case B, albeit with a significant rise in the velocity.

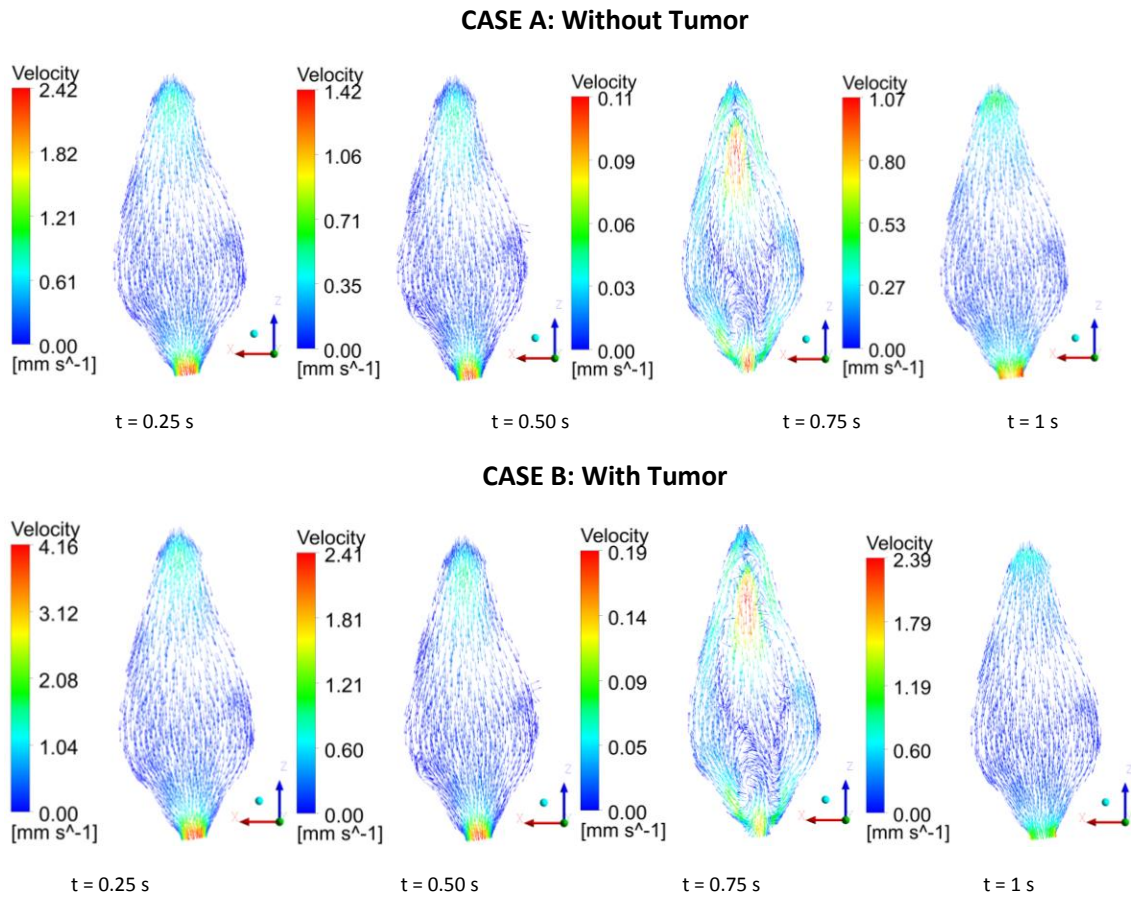

**Figure 8.** CSF flow in fourth ventricle during systole and diastole. Case A represents velocity vectors in fourth ventricles for a normal case where no tumor is present; and Case B represents velocity in fourth ventricle under the influence of tumor.

*3.5. Effect of Changing Flow of CSF on Ventricle Body—Deformation*

The deformation due to the interaction of fluid forces and tumor across the walls of ventricles are recorded. Figure 9 (a) and (b) shows a graph which gives deformation in time for both cases. For normal case the peak deformation in our study is found to be 2 um, whereas for the tumor case the peak deformation is found across third ventricle which came out as 11 um. This deformation is passed on to the CSF domain, where relative displacement of walls inwardly causes a rise in CSF flow velocity in third ventricle and in cerebral aqueduct. Table 5 summarizes the results in this regard.

**Table 5** Deformation in ventricular system.

| Specimen | Peak Deformation | | |
|---|---|---|---|
| | Without tumor (mm) | With tumor (mm) | Increase in Deformation (mm) |
| Datapoint 1 | 0.0021 | 0.011 | 0.0089 |
| Datapoint 2 | 0.0032 | 0.016 | 0.0128 |
| Datapoint 3 | 0.0055 | 0.029 | 0.0235 |

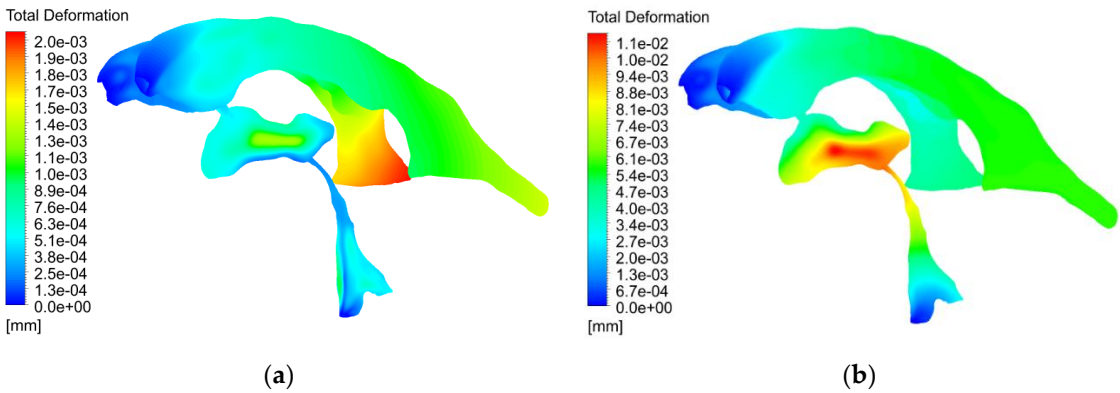

(**a**) (**b**)

Figure 9. Pressure on ventricle wall. (a) shows pressures without tumor; and (b) shows pressure under tumor forces

## 4. Discussion

This paper is one of the first attempts at modeling the CSF flow through the brain ventricles under the influence of brain tumor. By choosing to model the structure and fluid regimes using two-way Fluid-structure interaction, a complete understanding of the underlying dynamics of CSF flow in ventricles is presented. The proposed methodology considers all the relevant and major physiological factors which may markedly influence the CSF flow in the ventricles. This inter-alia includes cardiac cycle induced pulsations, dynamic interactions between the deformable wall of ventricles and flow profiles, backflow of CSF in diastole; and most importantly, the influence of brain tumor on physiological flow parameters of CSF. Among all the datapoints used in this paper, one clearly observes that under the case of brain tumor mean CSF velocity increases about 74.2%, which clearly shows the influence of brain tumor on CSF flow parametrization.

The results obtained from the proposed methodology are compared against the literature as shown in Table 6. Results are compared against two group of studies. First group corresponds to the numerical FSI studies done in past on the same connected matter and results herein are validated against these and second group corresponds PC-MRI studies which have also produced results regarding flow velocity of CSF et cetera. In all cases, peak CSF velocity is under the range reported by the proposed technique. For tumor, there are no established studies, hence we can resort this problem intuitively. Likewise, as evident from Table 6, our proposed model is not only compared against all the reported data in literature but also goes one step ahead. In studying CSF flow and biomechanical behavior, most papers in literature do not cater the deformable character of ventricles (excluding Masoumi et al [19]) and mostly rely on mere simple CFD-based flow modeling. However, our findings suggest that this assumption is not founded on firm roots and our proposed fully coupled FSI based model provides complete description of how ventricle and CSF interact dynamically. Second, based on this, we can easily infer that the proposed model can now (since it is fully two-way FSI) be used for various cases such as stenosis of cerebral aqueduct, effect of cerebral edema or hematoma on CSF flow et cetera.

Maximum local Reynold's number in the cerebral aqueduct reaches up to 107.91 which shows that the flow transition in the cerebral aqueduct is rather abrupt due to the narrow pathway. Furthermore, the proposition that the flow of CSF throughout the domain by and large remains a creeping flow is validated through this study as well. CSF flow velocity in lateral ventricles is of the orders of $10^{-3}$ mm/s. It is only through the aqueduct that the flow changes are maximum. Secondly, cardiac induced pulsations have an effect on the CSF flow profile mainly in the cerebral aqueduct and the fourth ventricle where backflow of cerebrospinal fluid occurs. Fourth ventricle usually experiences maximum backflow and hence during diastole velocity in craniocaudal direction reduces. This pattern repeats itself in all cardiac cycles for both tumor and non-tumor case; however, the only difference is the magnitude of velocity which is higher in case of tumor.

Transmantle pressures on the walls of ventricle remain under 4.3 Pascals (0.032 mmHg) for non-tumor case, whereas for tumor case, pressures increase up to 7.1 Pascals (0.053 mmHg) clearly highlighting that pressures rise during tumor case. The mean deformation of 0.0036 mm is recorded for the non-tumor cases, where no tumor forces are taken into consideration and only the fluid forces produce the aforementioned deformation. However, for the cases of tumor, mean deformation rise up to 0.0186 mm, signifying the fact that forces due to tumor constrict the ventricular wall and produce resultant increase in the CSF flow velocity. Needless to say, deformation depends upon the tumor volumes which are different in all cases, and so are the forces. But this nonetheless provides an insight of how much effect is produced by tumor on the deformation of ventricular walls. Lastly, the instant paper has dealt successfully in biomechanical modelling of CSF in human ventricular system under the influence of brain tumor. It first gives the foundational basis for modelling such scenarios using two-way FSI and then utilizing the modelling strategy it goes on to implement it for case of tumor and benchmarks it against both normal case and already reported established studies in literature. This biomechanical model of CSF in times of need for non-invasive techniques may be used by the doctors for the purpose of making informed decisions before going on for any surgical intervention.

Lastly, we would like to emphasize on the relevance of the results of this study. The magnitude of deformation obtained above may appear to be small, however, in biological tissues especially as delicate as brain is, this magnitude is significant. Not because the fact that tumor per-se is only the factor which itself produces such effect, but the fact that is a result of combined interaction. Consider a case when tumor is in close proximity to brainstem or cerebral aqueduct, such deformation can create havoc inside the aqueduct canal and may create increasing differentials of intracranial pressures resulting in tissue damage, hydrocephalus etc. In brain physiology, even a deformation of orders of millimeters is dangerous due to the delicate nature. Hence the results of this paper, when read with the overall neurology and physiological condition of brain provide the true extent of what results mean when one says CSF velocity has increased by a mean of 74.23%.

**Table 6** Comparison of results obtained from the proposed method against previous studies

| Study | Velocity in CA (mm/s) | Pressure drop (Pa) | Deformation (mm) |
|---|---|---|---|
| Jacobson et al. (1996) [12] | 28 | <1.1 for CA | Walls modeled as rigid |
| Fin and Grebe (2003) [13] | 64.65 | 2.91 for elastic wall model for CA | Implemented flexible cylindrical wall model of CA and obtained 0.61 mm deformation |
| Linninger et al. (2005) [15] | 25.8/-21.7 | 2 in CA | No deformation reported. However, net increase of ventricular volume of 4.5% is reported. |
| Kurtcuoglu et al. (2007)* [43] | 120 | 20 | No deformation reported |
| Howden et al. (2008) [18] | 11.38 | 1.14 | Walls modeled as rigid |
| Masoumi et al. (2010) [42] | 18/-15 | <5 | No deformation reported. However, net increase of ventricular volume of 6% is reported. |
| Masoumi et al. (2013) [19] | 8/-6 | < 2 | Normal CSF-Ventricular interaction under normal case: 0.006 mm |
| PC-CINE: Abbey et al (2009) [44] | 32.4±1.08 | - | - |
| PC-CINE: Algin et al (2010) [45] | 47.8±2.48 | - | - |
| PC-CINE: Lee et al (2004) [46] | 33.9±1.61 | - | - |
| **Proposed Method** | **Normal case: 22.8** | **Normal case: 2.88 Pa across CA** | **Normal case deformation= 0.002 mm** |
| | **Tumor case: 39.4** | **Tumor Case: 4.74 Pa across CA** | **Early tumor interaction deformation= 0.011 mm** |

**Author Contributions:** Muhammad Uzair Ul Haq and Ali Ahmed worked on the problem definition, methodology development, paper writing and Validation of Results, Faizan Nadeem provided validation expertise, Zartasha Mustansar and Arslan Shaukat acquired the funding and handled important topics related to paper definition and proof-reading, Asim Waris provided invaluable assistance in methodology development.

**Funding:** Please add: "This Research is funded by Higher Education Commission of Pakistan (HEC) under the auspices of National Research Program for Universities (NRPU) through an HEC grant code of 9954."

**Institutional Review Board Statement:** Not applicable

**Informed Consent Statement:** Not applicable

**Acknowledgments:** We would like to acknowledge INMOL cancer hospital, Lahore, Pakistan for providing MRI images for this study and National University of Sciences & Technology (NUST) for providing Computer Labs and logistic support necessary to conduct this study.

**Conflicts of Interest:** The authors declare no conflict of interest.